\documentclass[12pt]{amsart}
\usepackage{amsmath, amssymb, amsthm}
\usepackage{geometry}
\newtheorem{question}{Question}
\newtheorem{example}{Example}
\newtheorem{remark}{Remark}
\usepackage{enumitem}
\setlist[itemize]{leftmargin=2em}

\begin{document}

\title[An Affirmative Answer to Question 4 of Jachymski...]{An Affirmative Answer to Question 4 of Jachymski Concerning Mappings of Type $(\gamma,c)$}
\author{Maher Berzig}

\address{Maher Berzig\newline
\indent Universit\'e de Tunis \newline
\indent \'Ecole Nationale Sup\'erieure d'Ing\'enieurs de Tunis \newline
\indent D\'epartement de Math\'ematiques \newline 
\indent 5 Avenue Taha Hussein Montfleury, 1008 Tunis, Tunisia.}

\email{maher.berzig@gmail.com}

\date{}

\begin{abstract}
In his study of Cantor-type intersection theorems and fixed points of almost affine mappings, Jachymski \cite{Jachymski} posed four open questions. In this note, we provide an affirmative answer to the fourth question.
\end{abstract}



\maketitle
\section*{Introduction}

In his paper \cite{Jachymski}, Jachymski introduces the notion of mappings
of type $(\gamma,c)$ and poses four open questions. We focus on Question 4,
which asks:

\setcounter{question}{3}   
\begin{question}
Let $c\in(0,1)$. Does there exist a mapping of type $(\gamma,c)$, which is not
of type $(\gamma)$?
\end{question}

We recall the relevant definitions. Let
\[
\Gamma:=\{\gamma:\mathbb R_+\to\mathbb R_+:\gamma\text{ is strictly increasing,
continuous, convex, and }\gamma(0)=0\}.
\]
A mapping $T$ is of type $(\gamma,c)$ if for all $x,y$,
\[
\gamma\!\left(\left\| cTx+(1-c)Ty-T(cx+(1-c)y)\right\|\right)
\le
\|x-y\|-\|Tx-Ty\|.
\tag{1}
\]
A mapping is of type $(\gamma)$ if it is of type $(\gamma,c)$ for every
$c\in(0,1)$.

We prove that for every prescribed $c_0\in(0,1)$, there exists a mapping of
type $(\varphi,c_0)$ that is not of type $(\varphi)$. The construction uses
a piecewise-linear contraction with two different slopes and a linear
function $\varphi(s)=\lambda s\in\Gamma$, with $\lambda>1$ chosen to separate
the optimal constants for $c_0$ and a carefully chosen $c_*\ne c_0$.

The key observation is that the definition of $\Gamma$ imposes no upper bound
$\gamma(t)<t$; thus $\varphi(s)=\lambda s$ is valid for every $\lambda>0$.
This permits the separation argument that lies at the heart of the
construction.

\section{The construction}

We use the notation $\Gamma$ exactly as defined in \cite{Jachymski}:
\[
\Gamma:=\{\gamma\colon\mathbb R_+\to\mathbb R_+:\gamma\text{ is strictly increasing,
continuous, convex, and }\gamma(0)=0\}.
\]
In particular, no condition $\gamma(t)<t$ is imposed; thus $\varphi(s)=\lambda s$
with $\lambda>1$ is a valid element of $\Gamma$.

\subsection*{Construction of $T$}

Fix $c_0\in(0,1)$. Choose $c_*\in\left\{\frac13,\frac23\right\}$ with $c_*\ne c_0$, and set
\[
(a,b)=
\begin{cases}
\left(\frac15,\,\frac45\right), & c_*=\frac13,\\[2pt]
\left(\frac45,\,\frac15\right), & c_*=\frac23,
\end{cases}
\qquad
T(x)=\begin{cases} bx, & x<0,\\ ax, & x\ge 0.\end{cases}
\]
Since $0<a,b<1$, $T$ is a strict contraction (Lipschitz constant
$\max\{a,b\}<1$).

\subsection*{The Jensen defect and its supremum}

For $t\in(0,1)$ and $x\ne y$ set
\[
E_t(x,y):=tTx+(1-t)Ty-T(tx+(1-t)y),\qquad
D(x,y):=\|x-y\|-\|Tx-Ty\|.
\]
If $x,y$ have the same sign, $T$ is linear on their segment, so $E_t\equiv0$;
these pairs are irrelevant. If $x,y$ have opposite sign, write $x=-u<0<v=y$
($u,v>0$); an elementary computation splitting on the sign of $tx+(1-t)y$
gives
\[
E_t(x,y)=
\begin{cases}
t(a-b)u, & t\le v/(u+v),\\[2pt]
(1-t)(a-b)v, & t> v/(u+v),
\end{cases}
\qquad
D(x,y)=uA+vB>0,
\]
where $A:=1-b,\ B:=1-a$. Writing $r:=u/v$, the ratio $\|E_t\|/D$ equals
$t\|a-b\|\,r/(B+rA)$ for $r\le(1-t)/t$ and $(1-t)\|a-b\|/(B+rA)$ for
$r>(1-t)/t$; the first expression is increasing in $r$ and the second
decreasing, so the two meet — and the supremum over $r>0$ is attained — at
$r=(1-t)/t$. Substituting gives
\[
R(t):=\sup_{x\ne y}\frac{\|E_t(x,y)\|}{D(x,y)}
=\frac{\|a-b\|\,t(1-t)}{A(1-t)+Bt},
\]
attained exactly at the pair
\[
x=c_*-1,\quad y=c_*\qquad\text{when } t=c_*,
\]
i.e.\ $u=1-c_*,\ v=c_*$, so $u/v=(1-c_*)/c_*$ as required; one checks directly
that $c_*x+(1-c_*)y=c_*(c_*-1)+(1-c_*)c_*=0$.

\subsection*{$R$ has a unique maximum at $t=c_*$}

Differentiating $R$ and simplifying, the stationarity condition
$R'(t)=0$ reduces to $(B-A)t^2+2At-A=0$, whose positive root is
\[
t_*=\frac{\sqrt A}{\sqrt A+\sqrt B},\qquad\text{equivalently}\qquad
\frac{1-t_*}{t_*}=\sqrt{B/A}.
\]
For $c_*=\frac13$: $A=\frac15,\,B=\frac45$, so $B/A=4$ and $t_*=1/(1+2)=\frac13=c_*$. For
$c_*=\frac23$: $A=\frac45,\,B=\frac15$, so $B/A=\frac14$ and $t_*=1/(1+\frac12)=\frac23=c_*$. Since
$R(0^+)=R(1^-)=0$ and $R>0$ on $(0,1)$ with a single critical point, $t_*$ is
the unique maximum, so
\[
R(c_0)<R(c_*)\qquad\text{for } c_0\ne c_*.
\]

\subsection*{Conclusion of the construction}

Set $\lambda:=\tfrac12\big(1/R(c_*)+1/R(c_0)\big)$, so
$1/R(c_*)<\lambda<1/R(c_0)$, and let $\varphi(s):=\lambda s\in\Gamma$.
Since $D(x,y)>0$ for all $x\ne y$:
\begin{itemize}
\item For all $x\ne y$: $\varphi(\|E_{c_0}(x,y)\|)=\lambda\|E_{c_0}(x,y)\|\le
\lambda R(c_0)D(x,y)<D(x,y)$, since $\lambda R(c_0)<1$. Hence $T$ is of type
$(\varphi,c_0)$.
\item At the pair $x=c_*-1,\ y=c_*$ realizing $R(c_*)$: $\varphi(\|E_{c_*}
(x,y)\|)=\lambda R(c_*)D(x,y)>D(x,y)$, since $\lambda R(c_*)>1$. Hence the
defining inequality of type $(\varphi,c_*)$ fails for this pair, so $T$ is not
of type $(\varphi,c_*)$.
\end{itemize}

Since ``type $(\varphi)$'' means ``type $(\varphi,c)$ for every $c\in(0,1)$'',
and $T$ fails for $c=c_*$, $T$ is not of type $(\varphi)$. Hence the answer
to Question 4 is affirmative.

\begin{example}[Special case $c_0=\frac12$]
Take $c_*=\frac13$ (any choice with $c_*\ne \frac12$ works), so $(a,b)=\left(\frac15,\frac45\right)$:
\[
T(x)=\begin{cases}\dfrac{4}{5}x, & x<0,\\[7pt] \dfrac15 x, & x\ge0,\end{cases}
\]
and $\varphi(s)=\lambda s$ with $\lambda$ strictly between $1/R\left(\frac13\right)$ and
$1/R\left(\frac12\right)$. This gives an explicit mapping of type $(\varphi,\frac12)$ that is
not of type $(\varphi)$.
\end{example}

\begin{remark}
A key point is that $\lambda>1$ is permitted: the definition of $\Gamma$ in
\cite{Jachymski} imposes no condition $\gamma(t)<t$. Thus $\varphi(s)=\lambda s$
with $\lambda>1$ is a valid element of $\Gamma$. The construction therefore
works for every $c_0\in(0,1)$.
\end{remark}

\subsection*{Connection to Jachymski's Theorem 3.7}

Jachymski's Theorem 3.7 states that if $C$ is a nonempty closed bounded and
convex subset of a superreflexive Banach space $X$, and $T:C\to C$ is of type
$(\gamma,\frac12)$, then $T$ has a fixed point.
The proof of Theorem 3.7 relies crucially on the midpoint symmetry $c=\frac12$
to show that the approximate fixed-point sets satisfy the midpoint stability
condition of Proposition 3.3(iv). For this special value of $c$, the convex
combination $cx+(1-c)y$ becomes the midpoint $\frac{x+y}{2}$, which allows
the application of the almost-affine machinery developed in the paper.

Our construction in Question 4 shows that type $(\gamma,c_0)$ does
\emph{not} imply type $(\gamma)$ for any $c_0\in(0,1)$. This naturally
raises the following problem:

\begin{question}
Let $C$ be a nonempty closed bounded and convex subset of a superreflexive
Banach space $X$, and let $T:C\to C$ be of type $(\gamma,c)$ for a fixed
$c\in(0,1)$. Does $T$ necessarily have a fixed point?
\end{question}


\begin{thebibliography}{9}
\bibitem{Jachymski}
J. Jachymski,
``A Cantor type intersection theorem for superreflexive Banach spaces and
fixed points of almost affine mappings'',
\emph{J. Nonlinear Convex Anal.} \textbf{16} (2015), 1055--1068.
\end{thebibliography}
\end{document}